\documentclass[12pt,a4paper]{article}

\usepackage{latexsym}
\usepackage{a4wide}
\usepackage{amscd}
\usepackage{graphics}
\usepackage{amsmath}
\usepackage{amssymb}
\usepackage{soul,xcolor}
\usepackage{esint}
\usepackage{mathrsfs}
\usepackage{amsthm}
\usepackage{bbm}
\usepackage{color}
\usepackage{enumerate}
\usepackage{hyperref}
\usepackage[normalem]{ulem}
\input xy
\xyoption{all}

\newcommand{\N}{\mathbb{N}}                     % the natural numbers
\newcommand{\Z}{\mathbb{Z}}                     % the integer numbers
\newcommand{\R}{\mathbb{R}}                     % the real line
\newcommand{\C}{\mathbb{C}}                     % the complex plane
\newcommand{\T}{\mathbb{T}}                     % the torus
\newcommand{\im}{\mathrm{Im\,}}                 % immaginary part
\newcommand{\re}{\mathrm{Re\,}}                 % real part
\newtheorem{thm}{\sc Theorem}[section]               % numbered absolutely
\newtheorem*{thm*}{\sc Theorem}               % no number
\newtheorem*{cor*}{\sc Corollary}        % no number
\newtheorem{lem}[thm]{\sc Lemma}            % numbered along with Theorem
\newtheorem{prop}[thm]{\sc Proposition}     % numbered along with Theorem
\newtheorem{rem}[thm]{\sc Remark}           % numbered along with Theorem  
\title{A symplectic Hilbert hotel}

\author{Alberto Abbondandolo}

\date{}

\begin{document}

\maketitle

\begin{abstract}
We construct a smooth global symplectomorphism of an infinite dimensional separable Hilbert space mapping a ball of large radius to a subset of a cylinder of small radius. This symplectomorphism and its inverse map bounded sets to bounded sets. This example shows that Gromov's non-squeezing theorem does not generalize to infinite dimensional Hilbert spaces without further assumptions.
\end{abstract}

\section*{Introduction}

Regard $\C^n$ as a $2n$-dimensional real vector space with inner product $\langle \cdot, \cdot \rangle:= \re \langle \cdot, \cdot \rangle_{\C}$, norm $\|\cdot\|:=\langle \cdot ,\cdot \rangle^{\frac{1}{2}}$, and symplectic form $\omega_0:= \im \langle \cdot, \cdot \rangle_{\C}$, where $\langle z, w \rangle_{\C} = \sum_{j=1}^n z_j \bar{w}_j$ denotes the standard Hermitian product. A symplectomorphism between open subsets of $\C^n$ is a smooth diffeomorphism preserving $\omega_0$.
Denote by 
\[
B_r := \{ z\in \C^n \mid \|z\|<r \} \qquad \mbox{and} \qquad Z_r := \{ z=(z_1,\dots,z_n) \in \C^n \mid |z_1| < r \}  
\]
the open ball and open cylinder of radius $r$. Gromov's non-squeezing theorem \cite{gro85} asserts that if $R>r$ then no symplectomorphism can map $B_R$ to a subset of $Z_r$. 

Since the work of Kuksin \cite{kuk95,kuk95b} and Bourgain \cite{bou94}, there have been investigations on the question whether and to which extent the non-squeezing theorem extends to infinite dimensional Hilbert spaces. Recall that a (strong) symplectic form on a real Hilbert space $\mathbb{H}$ is a continuous skew-symmetric bilinear form $\omega_0: \mathbb{H} \times \mathbb{H} \rightarrow \R$ whose corresponding linear map $\mathbb{H} \rightarrow \mathbb{H}^*$ is an isomorphism. For the sake of concreteness, we here consider
\[
\mathbb{H} := \ell^2(\N,\C), \qquad \langle z,w \rangle_{\C} := \sum_{j\in \N} z_j \bar{w}_j,
\]
which we regard as a real Hilbert space with inner product $\langle \cdot, \cdot \rangle:= \re \langle \cdot, \cdot \rangle_{\C}$, norm $\|\cdot\|:=\langle \cdot ,\cdot \rangle^{\frac{1}{2}}$, and symplectic form $\omega_0:= \im \langle \cdot, \cdot \rangle_{\C}$. Every infinite dimensional separable real Hilbert space with a strong symplectic form is isomorphic to this model (see \cite[Section 2]{ama15}).  As in the finite dimensional setting, we denote by 
\[
B_r := \{ z\in \mathbb{H} \mid \|z\|<r \} \qquad \mbox{and} \qquad Z_r := \{ z=(z_1,\dots,z_n, \dots) \in \mathbb{H} \mid |z_1| < r \}  
\]
the open ball and open cylinder of radius $r$, and the question is whether symplectomorphisms of $(\mathbb{H},\omega_0)$ can squeeze $B_R$ into $Z_r$ when $R>r$. 

It is easy to show that non-squeezing holds for linear symplectic automorphisms of $(\mathbb{H},\omega_0)$, and  this remains true for nonlinear symplectomorphisms which are compact perturbations of linear ones \cite{kuk95,rig19b}. Remarkably, non-squeezing holds also for the time evolution of certain nonlinear Hamiltonian PDEs, even when their flow fails to be a compact perturbation of a linear flow. This is the case for the cubic Schr\"odinger equation on $\T:= \R/\Z$ \cite{bou94}, on $\R$ \cite{kvz19}, and on $\R^2$ \cite{kvz21}, the KdV equation on $\T$ \cite{ckstt05b} and on $\R$ \cite{nte22}, and the BBM equation on $\T$ \cite{rou10}. In all these papers, highly non-trivial finite dimensional reductions make it possible to apply Gromov's original non-squeezing theorem. There are also intrinsically infinite dimensional proofs of non-squeezing: Clarke's duality is used in \cite{ama15} to prove that non-squeezing holds when the image of the ball is convex, pseudoholomorphic curves are used in \cite{st16,st16b} to prove an abstract non-squeezing theorem with applications to the discrete nonlinear Schr\"odinger equation, and in \cite{cic19,cic21} to prove another interesting abstract non-squeezing theorem which recovers several of the previously mentioned results on Hamiltonian PDEs.

\medskip

The aim of this note is to construct a global smooth symplectomorphism of $\mathbb{H}$ that squeezes the ball of radius $R$ into the cylinder of radius $r$, where $0<r<R$ are fixed. More precisely, we prove the following result.

\begin{thm*}
Let $0<r<R$. There exists a smooth global symplectomorphism $\varphi: \mathbb{H} \rightarrow \mathbb{H}$ such that
\[
\varphi(B_R) \subset Z_r.
\]
Moreover, both $\varphi$ and $\varphi^{-1}$ map bounded sets to bounded sets.
\end{thm*}

Therefore, without further assumptions the non-squeezing theorem does not hold on infinite dimensional symplectic Hilbert spaces. 
The symplectomorphism $\varphi$ which we construct in the proof of the above theorem is Lipschitz continuous on $B_{\rho}$ for every $\rho<r$, but is not uniformly continuous on $B_r$ (see Remarks \ref{r:lipschitz} and \ref{r:non-unif-cont} below). In particular, its differential is not bounded on $B_r$. Symplectomorphisms arising from Hamiltonian PDEs usually have good derivative bounds. It is an interesting open question whether non-squeezing holds for symplectomorphisms of $\mathbb{H}$ all of whose derivatives are bounded on bounded sets.

\medskip

We conclude this introduction by describing some of the ideas involved in the construction of the example which proves  the above theorem. The starting idea is the first property of the Hilbert hotel \cite{hil13}: moving the guest in room $n$ to room $n+1$ for every $n\geq 1$ frees the first room. Said otherwise, the shift operator $(z_1,\dots, z_n, \dots)\mapsto (0,z_1,\dots, z_n, \dots)$, which preserves $\omega_0$, maps the whole $\mathbb{H}$ to a subset of $Z_r$, for every $r>0$. Of course, this map is not a genuine symplectomorphism because it is not open. However, it is not difficult to turn it into a symplectomorphism $\vartheta$ mapping the whole $\mathbb{H}$ to an open subset of the cylinder $Z_r$.

Here is the construction. We start by considering a norm-preserving symplectomorphism $\theta :\C^2 \rightarrow \C^2$ which is the identity near the origin and coincides with the  map $(z_1,z_2)\mapsto (z_2,-z_1)$ outside of a ball of radius $r$. For instance, $\theta$ could be the time-1 map of the Hamiltonian system defined by the Hamiltonian
\[
H(z_1,z_2) = \frac{\pi}{2} \chi(|z_1|^2+|z_2|^2) \,\mathrm{Im}\, ( \bar{z}_1 z_2 ) \qquad \forall (z_1,z_2)\in \C^2,
\]
where the smooth cut-off function $\chi$ satisfies $\chi(\rho)=0$ for $\rho$ close to 0 and $\chi(\rho)=1$ for $\rho\geq r^2$ (see Section \ref{s:psfc} below for sign conventions on Hamiltonian vector fields).  These properties imply that for every $(z_1,z_2)\in \C^2$ we have
\begin{equation}
\label{e:key}
|z_2|< r \qquad \Rightarrow \qquad |\pi_1 \theta(z_1,z_2) | < r,
\end{equation}
where $\pi_1: \C^2 \rightarrow \C$ denotes the projection onto the first factor. We denote by $\theta_n : \mathbb{H} \rightarrow \mathbb{H}$ the symplectomorphism which is obtained by acting by $\theta$ on the $n$-th and $(n+1)$-th coordinates. Then we define $\vartheta: \mathbb{H} \rightarrow \mathbb{H}$ as the limit
\[
\vartheta := \lim_{n\rightarrow \infty} \theta_1 \circ \theta_2 \circ \dots \circ \theta_n.
\]
It is easy to show that the above limit stabilizes locally, and hence defines a smooth symplectomorphism onto some open subset of $\mathbb{H}$. Property \eqref{e:key} implies that $\vartheta(\mathbb{H})$ is contained in $Z_r$. Indeed, any $z\in \mathbb{H}$ satisfies $|z_{n+1}|<r$ for $n$ large enough, so \eqref{e:key} implies that $|\pi_n \theta_n(z)|<r$, which implies $|\pi_{n-1} \theta_{n-1}\circ \theta_n(z)|<r$, and iterating this we obtain that 
\[
|\pi_1 \vartheta(z)| = |\pi_1 \theta_1 \circ \theta_2 \circ \dots \circ \theta_n(z) | < r.
\]

The drawback of this simple example is that it is not a global symplectomorphism of $\mathbb{H}$, as it is not surjective. Actually, it cannot be modified outside of a large ball $B_R$ with $R>r$ so that it becomes a global symplectomorphism of $\mathbb{H}$ mapping $B_R$ to a subset of $Z_r$. Indeed, if $\widetilde{\vartheta}: \mathbb{H} \rightarrow \mathbb{H}$ is a homeomorphism coinciding with $\vartheta$ on $B_R$, and hence on its closure, then the fact that homeomorphisms preserve boundaries and the fact that $\vartheta$ is norm-preserving imply
\[
\partial \vartheta(B_R) = \partial \widetilde{\vartheta}(B_R) = \widetilde{\vartheta} (\partial B_R) = \vartheta (\partial B_R) \subset \partial B_R,
\]
showing that the open set $\vartheta(B_R)$ is relatively closed in $B_R$. Since the latter space is connected, this implies that $\vartheta(B_R)=B_R$, contradicting $\vartheta(B_R)\subset Z_r$.

The global symplectomorphism $\varphi$ of the above theorem will be given by a composition $\varphi=\psi\circ \tau$, where $\tau$ translates the second coordinate $z_2$ by $R+s$, and $\psi$ is a norm-preserving global symplectomorphism of $\mathbb{H}$ mapping the complement of the cylinder $\{|z_2|<s\}$ to a subset of $Z_r$. Here, $s$ is any given positive number. The main point is to construct a map $\psi$ with the above properties. 

Here is where we make use of the second property of the Hilbert hotel: the hotel contains two copies of itself, by considering rooms with odd and even numbers. The construction of $\psi$ uses the splitting of $\mathbb{H}$ into the subspaces spanned by odd and even coordinates. The even coordinates act as control variables, damping the effect of the transformation on the odd coordinates. This is what makes surjectivity possible. For better bookkeeping, we shall work with $\ell^2(\Z,\C)$ instead of $\ell^2(\N,\C)$, and distinguish between positive and negative indices instead of odd and even ones.   

\paragraph{Acknowledgments.} This work was done while the author was in residence at the Institute for Advanced Study during a Special Year on Conformally Symplectic Dynamics and Geometry, supported by an anonymous donor. The AI tools provided by the IAS gave substantial assistance in finding the example constructed in this paper. The exposition is human generated. 

\section{Setting and further comments}
\label{s:psfc}

It is more convenient to work with the Hilbert space
\[
\mathbb{H} := \ell^2(\Z,\C),
\]
which is of course isomorphic to $\ell^2(\N,\C)$. The real inner product $\langle \cdot , \cdot \rangle$ and the symplectic form $\omega_0$ are the real and imaginary part of the standard Hermitian product 
\[
\langle  z,w \rangle_{\C} := \sum_{j\in \Z} z_j \bar{w}_j
\]
on $\mathbb{H}$.  We use the following convention for the vector field $X_H$ which is induced by a smooth Hamiltonian $H$ on $(\mathbb{H},\omega_0)$:
\[
\omega_0 (X_H,\cdot ) = dH.
\]
Since $\omega_0(x,y) = \langle -ix,y \rangle$ and $dH = \langle \nabla H,\cdot \rangle$, this convention reads
\begin{equation}
\label{e:ham}
X_H = i \nabla H = 2 i \frac{\partial H}{\partial \bar{z}},
\end{equation}
where the Wirtinger derivative $\frac{\partial}{\partial \bar{z}}$ acts componentwise. When $X_H$ is complete, we denote by $\phi^t_H$ its flow. The Poisson bracket of two Hamiltonians reads
\begin{equation}
\label{e:poisson}
\{H,K\} := \omega_0 (X_H,X_K) = 4 \, \mathrm{Im}\, \Bigl\langle \frac{\partial H}{\partial \bar{z}} , \frac{\partial K}{\partial \bar{z}} \Bigr\rangle_{\C}.
\end{equation}
As cylinder $Z_r$ we now consider the set
\[
Z_r := \{ z \in \mathbb{H} \mid |z_0|< r\}.
\]
Given $s>0$, we also consider the set
\[
C_s := \{ z\in \mathbb{H} \mid |z_{-1}| \geq s\}.
\]
As anticipated in the Introduction, the construction of the squeezing symplectomorphism is based on the following result.

\begin{prop}
\label{p:main}
For every $r>0$ and $s>0$ there exists a smooth norm-preserving global symplectomorphism $\psi: \mathbb{H} \rightarrow \mathbb{H}$ such that
\[
\psi(C_s) \subset Z_r.
\]
\end{prop}

The theorem stated in the introduction is an immediate consequence of the above proposition. Indeed, let $\tau: \mathbb{H} \rightarrow \mathbb{H}$ be the translation adding $R+s$ to the component $z_{-1}$, and set
\[
\varphi:= \psi\circ \tau.
\]
This map is a global symplectomorphism. Since $\tau$ maps $B_R$ to a subset of $C_s$, Proposition \ref{p:main} implies that $\varphi(B_R)\subset Z_r$.  Since $\tau$ and $\tau^{-1}$ are isometries and since $\psi$ and $\psi^{-1}$ are norm preserving, both $\varphi$ and $\varphi^{-1}$ map bounded sets to bounded sets. After identifying $\ell^2(\Z,\C)$ with $\ell^2(\N,\C)$ by a linear symplectic isometry identifying the coordinate of index $0\in \Z$ with the coordinate of index $1\in \N$, this concludes the proof of the theorem.

\begin{rem}
{\rm In this example, we also have
\[
\varphi(B_R) \subset B_{2R+s},
\]
where $s$ can be arbitrarily small.}
\end{rem}
 
\begin{rem}
\label{r:Hamiltonian}
{\rm The symplectomorphism $\varphi$ is Hamiltonian, meaning that there exists a smooth time-dependent Hamiltonian function $H:[0,1]\times \mathbb{H} \rightarrow \R$ such that the time-one map of the corresponding time-dependent Hamiltonian vector field is $\varphi$. Indeed, this is true for every global symplectomorphism of $(\mathbb{H},\omega_0)$, as the usual finite-dimensional proof based on Alexander's trick works also on infinite dimensional Hilbert spaces. We sketch the argument for the reader's convenience. The formula
\[
\psi_t(z) := \frac{\varphi(tz)-\varphi(0)}{t} = \int_0^1 d\varphi_{stz}(z) \, ds
\]
defines a smooth path of global symplectomorphisms of $(\mathbb{H},\omega_0)$ such that $\psi_1 = \varphi - \varphi(0)$ and $\psi_0=d\varphi_0$ is linear. Recall that the linear symplectic group of $(\mathbb{H},\omega_0)$ is connected, because the polar decomposition gives us a homotopy equivalence with the unitary group of $\mathbb{H}$ seen as complex Hilbert space, which is connected (and even contractible, by Kuiper's theorem). Therefore, from $\{\psi_t\}_{t\in [0,1]}$ we easily obtain a smooth path of symplectomorphisms $\{\varphi_t\}_{t\in [0,1]}$ such that $\varphi_0$ is the identity and $\varphi_1=\varphi$. The fact that $\varphi_t$ preserves $\omega_0$ implies that the contraction of $\omega_0$ by the vector field $X_t := (\frac{d}{dt} \varphi_t ) \circ \varphi_t^{-1}$ is a closed one-form. Since $\mathbb{H}$ is simply connected, this one-form is exact and we obtain a smooth one-parameter family of Hamiltonians $\{H_t\}_{t\in [0,1]}$ such that $\omega_0(X_t,\cdot)=dH_t$.}
\end{rem}

\begin{rem}
\label{r:autonomous}
{\rm It is also possible to squeeze the ball of radius $R$ into the cylinder of radius $r$ by the time-one map of an autonomous Hamiltonian vector field. This can be achieved by a lifting procedure to the time-energy extended phase space. Indeed, assume that $\varphi: \mathbb{H} \rightarrow \mathbb{H}$ is a global symplectomorphism mapping  $B_R$ to a bounded subset of $Z_r$.  As seen in Remark \ref{r:Hamiltonian}, there is a smooth time-dependent Hamiltonian $H: [0,1]\times \mathbb{H} \rightarrow \R$ whose non-autonomous flow $\varphi_t$ satisfies $\varphi_1 = \varphi$. Consider the extended phase space $\widehat{\mathbb{H}} := \R^2 \oplus \mathbb{H}$ with coordinates $(x,y,z)$ and symplectic form $\widehat{\omega}_0 = (dx\wedge dy) \oplus \omega_0$. The space $(\widehat{\mathbb{H}},\widehat{\omega}_0)$ is symplectically and isometrically equivalent to $(\mathbb{H},\omega_0)$, and it suffices to construct an autonomous Hamiltonian vector field on $\widehat{\mathbb{H}}$ whose time-one map squeezes the ball of radius $R$ into a cylinder of radius $r$. On $\widehat{\mathbb{H}}$, we consider the autonomous Hamiltonian
\[
\widehat{H}(x,y,z) := 3R \bigl( y + \rho'(x) H(\rho(x),z) \bigr),
\]
where $\rho: \R \rightarrow [0,1]$ is a smooth function having the value $0$ on $(-\infty,R]$ and the value $1$ on $[2R,+\infty)$. The curve $t\mapsto (x(t),y(t),z(t))$ is an integral line of $X_{\widehat{H}}$ if and only if
\begin{eqnarray}
\label{e:s1}
\dot{x} &=& \frac{\partial \widehat{H}}{\partial y}(x,y,z) = 3R, \\ 
\label{e:s2}
\dot{y} &=& - \frac{\partial \widehat{H}}{\partial x}(x,y,z) = - 3 R \Bigl(\rho''(x) H(\rho(x),z) + \rho'(x)^2 \frac{\partial H}{\partial t} (\rho(x),z) \Bigr), \\
\label{e:s3}
\dot{z} &=& 3 R \rho'(x) X_{H({\rho(x)},\cdot)} (z).
\end{eqnarray}
The solution of (\ref{e:s1}--\ref{e:s3}) with $x(0)=x_0$, $y(0)=y_0$, and $z(0)=z_0$ is
\begin{eqnarray}
\label{e:sol1}
x(t) &=& x_0 + 3Rt, \\ 
\label{e:sol2}
y(t) &=& y_0 + \rho'(x_0) H(\rho(x_0),z_0) - \rho'(x(t)) H(\rho(x(t)),z(t)), \\
\label{e:sol3}
z(t) &=& \varphi_{\rho(x(t))} \circ \varphi_{\rho(x_0)}^{-1} (z_0).
\end{eqnarray}
Indeed, \eqref{e:sol1} and \eqref{e:sol3} are obtained by integrating \eqref{e:s1} and \eqref{e:s3}, and \eqref{e:sol2} follows from the energy conservation 
\[
\widehat{H}(x(t),y(t),z(t)) = \widehat{H}(x_0,y_0,z_0).
\]
This shows that $X_{\widehat{H}}$ is complete. If $(x_0,y_0,z_0)$ belongs to the ball $\widehat{B}_R$ of radius $R$ in $\widehat{\mathbb{H}}$ then $|x_0|<R$ and hence $x(1)>2R$, implying that $\rho'(x_0)=\rho'(x(1))=0$. Therefore, (\ref{e:sol1}--\ref{e:sol3}) give 
\[
\widehat{\varphi} (x_0,y_0,z_0) := \phi^1_{\widehat{H}}(x_0,y_0,z_0) = \bigl(x_0+3R,y_0,\varphi(z_0)\bigr).
\]
The inclusion $\varphi(B_R)\subset Z_r$ then gives us $\widehat{\varphi}(\widehat{B}_R)\subset \widehat{Z}_r$, where $\widehat{Z}_r:= \R^2 \times Z_r$. Therefore, the autonomous Hamiltonian diffeomorphism $\widehat{\varphi}$ squeezes the ball of radius $R$ into a cylinder of radius $r$. Moreover, $\widehat{\varphi}(\widehat{B}_R)$ is bounded.}
\end{rem}

\section{Building blocks}

The remaining part of this note is devoted to the proof of Proposition \ref{p:main}. The numbers $r>0$ and $s>0$ are the constants fixed in its statement. On $\mathbb{H}$, we consider the quadratic Hamiltonian
\[
K_n(z) := \frac{\pi}{2} \mathrm{Im}\, \bigl( \bar{z}_0 z_n + \bar{z}_{-n} z_{-n-1} \bigr) = \frac{\pi}{4i} ( \bar{z}_0 z_n - z_0 \bar{z}_n + \bar{z}_{-n} z_{-n-1} - z_{-n} \bar{z}_{-n-1} ).
\]
By \eqref{e:ham}, its Hamiltonian vector field has the form
\[
X_{K_n} (z) = \frac{\pi}{2} J_n z,
\]
where $J_n: \mathbb{H} \rightarrow \mathbb{H}$ is the skew-adjoint operator mapping the coordinates $(z_{-n-1},z_{-n},z_0,z_n)$ to $(-z_{-n}, z_{-n-1}, z_n, - z_0)$ and the other coordinates to zero. The flow of $X_{K_n}$ is the unitary flow
\begin{equation}
\label{e:flowK}
\phi_{K_n}^t(z) = e^{ \frac{\pi}{2} t J_n} z, \qquad \forall (t,z)\in \R \times \mathbb{H}.
\end{equation}
In particular,
\begin{equation}
\label{e:mapK}
\phi^1_{K_n} (z) = U_n z,
\end{equation}
where $U_n:= e^{\frac{\pi}{2} J_n}$ is the unitary operator mapping the coordinates $(z_{-n-1},z_{-n},z_0,z_n)$ to $(-z_{-n}, z_{-n-1}, z_n, - z_0)$ and leaving the other coordinates unchanged.

Let $(a_n)_{n\geq 1}$  be a sequence such that
\[
0<a_n < 1, \qquad (a_n) \mbox{ is nondecreasing}, \qquad a_n\rightarrow 1.
\] 
For every integer $n\geq 1$, let $\chi_n: \R \rightarrow [0,1]$ be a smooth function such that
\[
\chi_n = 0 \quad \mbox{on } (-\infty,a_n] , \qquad \chi_n=1 \quad \mbox{on } [1,+\infty).
\] 
On $\mathbb{H}$, we consider the quadratic functions
\[
p_n(z):= \frac{1}{r^2} \Bigl( |z_0|^2 + \sum_{j\geq n} |z_j|^2 \Bigr), \qquad q_n(z) := \frac{1}{s^2} \sum_{j\leq -n} |z_j|^2, 
\]
and the smooth Hamiltonian 
\[
H_n := \chi_n(p_n) \chi_n (q_n) K_n.
\]

\begin{lem}
\label{l:psin}
The flow of $X_{H_n}$ is complete and preserves the functions $K_n$, $p_n$, $q_n$,  and the norm of $\mathbb{H}$. The symplectomorphism
\begin{equation}
\label{e:psin}
\psi_n := \phi_{H_n}^1
\end{equation}
satisfies:
\begin{enumerate}[(a)]
\item if $p_n(z)\geq 1$ and $q_n(z)\geq 1$ then $\psi_n(z) = U_n z$;
\item if $p_n(z)\leq a_n$ or $q_n(z)\leq a_n$ then  $\psi_n(z)=z$;
\item $\psi_n$ changes the coordinates $(z_{-n-1},z_{-n},z_0,z_n)$, multiplies the coordinates $z_j$ with $j\leq -n-2$ or $j\geq n+1$ by unimodular factors, and leaves the coordinates $z_j$ with $1\leq |j| \leq n-1$ unchanged.
\end{enumerate}
\end{lem}

\begin{proof}
By \eqref{e:ham}, the Hamiltonian vector fields of the quadratic functions $p_n$ and $q_n$ are
\[
X_{p_n} (z) = \frac{2i}{r^2} \Pi_n^+ z, \qquad X_{q_n}(z) = \frac{2i}{s^2} \Pi_n^- z \qquad \forall z\in \mathbb{H},
\]
where $\Pi_n^+$ is the orthogonal projector onto the subspace corresponding to the coordinates $z_j$ with $j=0$ or $j\geq n$, while $\Pi_n^-$ is the orthogonal projector onto the subspace corresponding to the coordinates $z_j$ with $j\leq -n$. Their flows are the unitary flows
\[
\phi_{p_n}^t(z) = (\mathrm{id} - \Pi_n^+) z + e^{\frac{2it}{r^2}} \Pi_n^+ z, \qquad
\phi_{q_n}^t(z) = (\mathrm{id} - \Pi_n^-) z + e^{\frac{2it}{s^2}} \Pi_n^- z.
\]
The first flow multiplies the coordinates $z_j$ with $j=0$ or $j\geq n$ by the unimodular factor $e^{\frac{2it}{r^2}}$ and leaves the other coordinates invariant. The second flow multiplies the coordinates $z_j$ with $j\leq -n$ by the unimodular factor $e^{\frac{2it}{s^2}}$ and leaves the other coordinates invariant. These two flows clearly commute. They also commute with the Hamiltonian flow of $K_n$ given by \eqref{e:flowK} because by \eqref{e:poisson} we have
\[
\begin{split}
\{K_n,p_n\}  &= 4 \, \mathrm{Im}\, \Bigl\langle \frac{\partial K_n}{\partial \bar{z}} , \frac{\partial p_n}{\partial \bar{z}} \Bigr\rangle_{\C}= \frac{\pi}{r^2} \, \mathrm{Im}\, \Bigl( \frac{z_n \bar{z}_0 - z_0 \bar{z}_n}{i} \Bigr)=0, \\
\{K_n,q_n\}  &= 4 \, \mathrm{Im}\, \Bigl\langle \frac{\partial K_n}{\partial \bar{z}} , \frac{\partial q_n}{\partial \bar{z}} \Bigr\rangle_{\C}= \frac{\pi}{s^2} \, \mathrm{Im}\, \Bigl( \frac{z_{-n-1} \bar{z}_{-n} - z_{-n} \bar{z}_{-n-1}}{i} \Bigr)=0.
\end{split}
\]
Since moreover $\{p_n,q_n\}=0$, all Poisson brackets among $K_n$, $p_n$, and $q_n$ vanish, and hence each of the three linear flows preserves the three functions $K_n$, $p_n$, and $q_n$. The Hamiltonian vector field of $H_n$ has the form
\[
X_{H_n} = \chi_n(p_n) \chi_n(q_n) X_{K_n} + \chi_n'(p_n) \chi_n(q_n) K_n X_{p_n} + \chi_n(p_n) \chi_n'(q_n) K_n X_{q_n}.
\]
Since $K_n$, $p_n$, and $q_n$ are constant along the three linear flows, the coefficients in the above expression are constant along the integral curves of $X_{H_n}$. Together with the commutativity of the three linear flows, this implies that $X_{H_n}$ is complete and its flow has the form
\begin{equation}
\label{e:flowH}
\phi_{H_n}^t = \phi_{K_n}^{\chi_n(p_n)\chi_n(q_n)t} \circ \phi_{p_n}^{\chi_n'(p_n) \chi_n(q_n) K_n t} \circ \phi_{q_n}^{ \chi_n(p_n) \chi_n'(q_n) K_n t}.
\end{equation}
This flow preserves $K_n$, $p_n$, $q_n$,  and also the norm, because the three linear flows are unitary. If $z\in \mathbb{H}$ is such that $p_n(z)\geq 1$ and $q_n(z)\geq 1$ then
\[
\chi_n(p_n(z))=1, \quad \chi_n'(p_n(z))=0, \quad \chi_n(q_n(z))=1, \quad \chi_n'(q_n(z))=0,
\]
and hence
\[
\psi_n(z) = \phi^1_{H_n}(z) = \phi_{K_n}^1(z) = U_n z
\]
by \eqref{e:mapK}, proving (a). If $p_n(z)\leq a_n$ or $q_n(z)\leq a_n$ then 
\[
\chi_n(p_n(z))\chi_n(q_n(z))=\chi_n'(p_n(z))\chi_n(q_n(z)) = \chi_n(p_n(z))\chi_n'(q_n(z))=0,
\] 
and hence
\[
\psi_n(z) = z,
\]
proving (b). Statement (c) follows from \eqref{e:flowH} and the properties of the three linear flows discussed above.
\end{proof}

\section{Construction of the symplectomorphism $\psi$}

Set $\Psi_0:= \mathrm{id}$ and
\[
\Psi_n := \psi_n \circ \dots \circ \psi_1,
\]
where $\psi_k$ is the smooth symplectomorphism defined in \eqref{e:psin}.Our first aim is to show that this sequence stabilizes locally and its limit for $n\rightarrow \infty$  defines a smooth global symplectomorphism of $\mathbb{H}$. 

For every $n\geq 1$, we consider the functions
\[
P_n(z) := r^2 p_n(\Psi_{n-1}(z))  \qquad Q_n(z) := s^2 q_n( \Psi_{n-1}(z)).
\]
By the definition of $p_{n+1}$, we have
\begin{equation}
\label{e:bound}
|\pi_0 \Psi_n(z)|^2 \leq P_{n+1}(z).
\end{equation}
Here and in what follows, $\pi_k : \mathbb{H} \rightarrow \C$ denotes the projection to the $k$-th coordinate, for $k\in \Z$. 
The next two Lemmata establish some useful properties of these functions. 

\begin{lem}
\label{l:monotone}
For every $z\in \mathbb{H}$ and $n\geq 1$ we have
\[
P_{n+1}(z) = P_n(z) - |\pi_n \Psi_n(z)|^2, \qquad  Q_{n+1}(z) = Q_n(z) - |\pi_{-n} \Psi_n(z)|^2.
\]
In particular, the sequences $(P_n(z))_{n\geq 1}$ and $(Q_n(z))_{n\geq 1}$ are nonincreasing.
\end{lem}

\begin{proof}
From the identity $p_{n+1}(z) = p_n(z) - \frac{1}{r^2} |z_n|^2$ and from the fact that $\psi_n$ preserves the function $p_n$, by Lemma \ref{l:psin}, we deduce 
\[
\begin{split}
P_{n+1}(z) &= r^2 p_{n+1}(\Psi_n(z)) = r^2 p_n (\psi_n \circ \Psi_{n-1}(z)) - |\pi_n \Psi_n(z)|^2 \\ &= r^2 p_n(\Psi_{n-1}(z)) - |\pi_n \Psi_n(z)|^2 = P_n(z) - |\pi_n \Psi_n(z)|^2,
\end{split}
\]
proving the first identity. The second identity is proven similarly, starting from $q_{n+1}(z) = q_n(z) - \frac{1}{s^2} |z_{-n}|^2$ and using the fact that $\psi_n$ preserves $q_n$.  
\end{proof}

\begin{lem}
\label{l:lin}
Let $z\in \mathbb{H}$ be such that $\Psi_n(z) = U_n \Psi_{n-1}(z)$ for every $n\geq 1$. Then
\[
P_{n+1}(z) = \sum_{j\geq n} |z_j|^2 \qquad \forall n\geq 1,
\]
and in particular $P_n(z)\rightarrow 0$.
\end{lem}

\begin{proof}
The operator $U_k$ changes only the coordinates $z_{-k-1}, z_{-k}, z_0, z_k$. By assumption, $\Psi_n(z) = U_n \circ \dots \circ U_1 z$ and hence the coordinates of index $j\geq n+1$ of $\Psi_n(z)$ are unchanged, while 
\[
\pi_0 \Psi_n(z) = \pi_0 U_n \circ \dots \circ U_1 z  = \pi_n U_{n-1} \circ \dots \circ U_1 z = z_n.
\]
Therefore,
\[
P_{n+1}(z) = |\pi_0 \Psi_n(z)|^2 + \sum_{j\geq n+1} |\pi_j  \Psi_n(z)|^2 = \sum_{j\geq n} |z_j|^2.
\]
\end{proof}

In the next lemma, we show that the sequences $\Psi_n$ and $\Psi_n^{-1}$ stabilize locally.

\begin{lem}
\label{l:stabilizes}
For every $z\in \mathbb{H}$ there are a neighborhood $U$ of $z$ and $N\in \N$ such that
\[
\Psi_n|_U = \Psi_N|_U \qquad \mbox{and} \qquad \Psi_n^{-1}|_U = \Psi_N^{-1}|_U
\]
for every $n\geq N$.
\end{lem}

\begin{proof}
We claim that there is $n_0\geq 1$ such that $P_{n_0}(z)< a_{n_0} r^2$ or $Q_{n_0}(z) < a_{n_0} s^2$. 

Suppose for contradiction that $P_n(z)\geq  a_n r^2$ and $Q_n(z) \geq a_n s^2$ for every $n\geq 1$. Since the sequences $(P_n(z))$ and $(Q_n(z))$ are nonincreasing and $(a_n)$ converges to 1, we deduce that
\[
P_n(z) \geq r^2 \qquad \mbox{and} \qquad Q_n(z) \geq s^2 \qquad \forall n\geq 1.
\]
Therefore,
\[
p_n(\Psi_{n-1}(z)) \geq 1  \qquad \mbox{and} \qquad q_n(\Psi_{n-1}(z)) \geq 1  \qquad \forall n\geq 1,
\]
and hence Lemma \ref{l:psin} (a) implies that $\Psi_n(z) = U_n \Psi_{n-1}(z)$ for every $n\geq 1$. By Lemma \ref{l:lin}, $(P_n(z))$ tends to zero, which is a contradiction. This proves the claim.

We now assume that there is $n_0\geq 1$ such that $P_{n_0}(z)< a_{n_0} r^2$, the case $Q_{n_0}(z) < a_{n_0} s^2$ being analogous. Since the function $P_{n_0}$ is continuous, we find a neighborhood $U_1$ of $z$ on which $P_{n_0} < a_{n_0} r^2$. For $w\in U_1$ and $n\geq n_0$, Lemma \ref{l:monotone} and the monotonicity of $(a_n)$ give
\[
r^2 p_n(\Psi_{n-1}(w)) = P_n(w) \leq P_{n_0}(w) < a_{n_0} r^2  \leq a_n r^2.
\]
Therefore, $p_n(\Psi_{n-1}(w)) < a_n$ and Lemma \ref{l:psin} (b) implies 
\[
\Psi_n(w) = \psi_n(\Psi_{n-1}(w)) = \Psi_{n-1}(w).
\]
This shows that $\Psi_n = \Psi_{n_0-1}$  on $U_1$ for every $n\geq n_0-1$.

Now consider the maps
\[
\Psi_n^{-1} = \psi_1^{-1} \circ \dots \circ \psi_n^{-1}.
\]
Choose $n_1\in \N$ such that 
\[
\sum_{j\leq -n_1} |z_j|^2 < a_1 s^2.
\]
If $U_2$ is a sufficiently small neighborhood of $z$, then
\[
\sum_{j\leq -n_1} |w_j|^2 < a_1 s^2 \qquad \forall w\in U_2.
\]
For $w\in U_2$ and $n\geq n_1$ we have
\[
q_n(w) = \frac{1}{s^2} \sum_{j\leq -n} |w_j|^2 \leq \frac{1}{s^2} \sum_{j\leq -n_1} |w_j|^2 < a_1 \leq a_n,
\]
and hence Lemma \ref{l:psin} (b) implies that $\psi_n(w)=w$, or equivalently $\psi_n^{-1}(w)=w$. This shows that $\psi_n^{-1}$ is the identity on $U_2$ for every $n\geq n_1$ and hence
\[
\Psi_n^{-1}|_{U_2} = \Psi_{n_1-1}^{-1}|_{U_2} \qquad \forall n\geq n_1.
\]
We conclude that the thesis holds with $U=U_1\cap U_2$ and $N=\max\{n_0,n_1\}-1$. 
  \end{proof}

The fact that $\Psi_n$ stabilizes locally implies that
\[
\psi := \lim_{n\rightarrow \infty} \Psi_n
\]
is a well-defined smooth map from $\mathbb{H}$ to $\mathbb{H}$ satisfying $\psi^* \omega_0 = \omega_0$. Since also $\Psi_n^{-1}$ stabilizes locally, the same is true for 
\[
\xi:= \lim_{n\rightarrow \infty} \Psi_n^{-1}.
\]
For every $z\in \mathbb{H}$, Lemma \ref{l:stabilizes} gives us natural numbers $N$ and $M$ such that $\psi(z) = \Psi_n(z)$ for every $n\geq N$ and $\xi(\psi(z)) = \Psi_m^{-1}(\psi(z))$ for every $m\geq M$. Choosing $n\geq \max\{N,M\}$, we find
\[
\xi\circ \psi(z) = \Psi_n^{-1} (\psi(z)) = \Psi_n^{-1} (\Psi_n(z)) = z,
\]
proving that $\xi\circ \psi=\mathrm{id}$ on $\mathbb{H}$. Similarly, $\psi\circ \xi=\mathrm{id}$ on $\mathbb{H}$, and hence $\psi$ is a global symplectomorphism with inverse $\xi$. Since all symplectomorphisms $\psi_n$ are norm preserving, by Lemma \ref{l:psin}, so is $\psi$.

\medskip

We now prove that
\begin{equation}
\label{e:main}
\psi(C_s) \subset Z_r.
\end{equation}
Let $z\in C_s$. By Lemma \ref{l:monotone}, the nonnegative sequence $(P_n(z))$ is nonincreasing and hence converges to some number $P_{\infty}\geq 0$. By Lemma \ref{l:stabilizes}, we have
\[
\psi(z) = \Psi_n(z)
\]
for every large integer $n$, and for such $n$ the bound \eqref{e:bound} reads
\[
|\pi_0 \psi(z)|^2 \leq P_{n+1}(z).
\]
By taking a limit in $n$, we deduce that
\[
|\pi_0 \psi(z)|^2 \leq P_{\infty},
\]
and hence it is enough to show that $P_{\infty} < r^2$.

Suppose, by contradiction, that $P_{\infty}\geq r^2$. Then $P_n(z)\geq r^2$ for every $n\geq 1$, and hence
\begin{equation}
\label{e:p}
p_n(\Psi_{n-1}(z)) = \frac{1}{r^2} P_n(z) \geq 1  \qquad \forall n\geq 1.
\end{equation}
We claim that, for every $n\geq 1$,
\begin{equation}
\label{e:claim}
| \pi_{-n} \Psi_{n-1}(z)| = |z_{-1}| \qquad \mbox{and} \qquad \Psi_n(z) = U_n \Psi_{n-1}(z).
\end{equation}
Note that the first identity in \eqref{e:claim} holds trivially for $n=1$. If this identity holds for some $n\geq 1$, then 
\[
q_n(\Psi_{n-1}(z)) \geq \frac{1}{s^2} | \pi_{-n} \Psi_{n-1}(z) |^2 = \frac{1}{s^2} |z_{-1}|^2 \geq  1 .
\]
By Lemma \ref{l:psin} (a), the above inequality and \eqref{e:p} imply that $\psi_n(\Psi_{n-1}(z)) = U_n  \Psi_{n-1}(z)$, which is the second identity in \eqref{e:claim}. Since $U_n$ swaps up to sign the coordinates with index $-n-1$ and $-n$, we now have
\[
|\pi_{-n-1} \Psi_n(z)| = |\pi_{-n-1} U_n \Psi_{n-1}(z)|  = |\pi_{-n} \Psi_{n-1}(z)| = |z_{-1}|,
\]
which is the first identity in \eqref{e:claim} for $n+1$. This proves \eqref{e:claim} for every $n\geq 1$.

Thanks to the second identity in \eqref{e:claim}, Lemma \ref{l:lin} gives $P_n(z)\rightarrow 0$, contradicting $P_{\infty}\geq r^2>0$. This proves \eqref{e:main}, concluding the proof of Proposition \ref{p:main}.

\medskip

In the next two remarks, $\{e_k\}_{k\in \Z}$ denotes the standard Hilbert basis of $\mathbb{H}$.

\begin{rem}
\label{r:lipschitz}
{\rm For every positive number $\rho <r$, the map $\psi$ is Lipschitz continuous on $\tau(B_{\rho})$, where $\tau(z)=z+ (R+s)e_{-1}$, and hence $\varphi=\psi\circ \tau$ is Lipschitz continuous on $B_{\rho}$, as claimed in the Introduction. Indeed, if $z=w+(R+s)e_{-1}$ with $w\in B_{\rho}$, then Lemma \ref{l:monotone} gives 
\[
P_n(z) \leq P_1(z) = \sum_{j\geq 0} |w_j|^2 \leq \|w\|^2 < \rho^2 \qquad \forall n\in \N.
\]
If $N\in \N$ is such that $a_N r^2 \geq \rho^2$, Lemma \ref{l:psin} (b) implies that $\psi=\Psi_{N-1}$ on $\tau(B_{\rho})$. Each norm-preserving map $\psi_k$ is a composition of unitary flows whose times are smooth functions of the quadratic forms $K_k$, $p_k$, and $q_k$, so $\psi_k$ has bounded differential on bounded sets. Hence the differential of $\Psi_{N-1}$ is bounded on $\tau(B_{\rho})$, and hence $\psi$ is Lipschitz continuous on $\tau(B_{\rho})$.}  
\end{rem}

\begin{rem}
\label{r:non-unif-cont}
{\rm The map $\psi$ is not uniformly continuous on $\tau(B_r)$, and hence $\varphi=\psi\circ \tau$ is not uniformly continuous on $B_r$, as claimed in the Introduction. Indeed, since uniform continuity passes to closures it is enough to show that $\psi$ is not uniformly continuous on $\tau(\overline{B}_r)$. The sequences
\[
w_n := r \sqrt{a_{n+1}} e_n + (R+s) e_{-1}, \qquad z_n:= r e_n + (R+s)e_{-1}
\]
take values in $\tau(\overline{B}_r)$. We claim that
\begin{equation}
\label{e:conto}
\pi_{n+1} \psi(w_n) = 0, \qquad  \pi_{n+1} \psi(z_n) = -r.
\end{equation}
The above identities imply that $\psi$ is not uniformly continuous on $\tau(\overline{B}_r)$. Indeed, $\|w_n-z_n\|= r (1-\sqrt{a_{n+1}})\rightarrow 0$ whereas \eqref{e:conto} implies that $\|\psi(w_n) - \psi(z_n)\|\geq r$. There remains to prove \eqref{e:conto}. From
\[
P_k(w_n)\leq P_1(w_n) = a_{n+1} r^2 \leq a_k r^2 \qquad \forall k\geq n+1
\]
we deduce by Lemma \ref{l:psin} (b) that $\Psi_k(w_n)=\Psi_n(w_n)$ for every $k\geq n$. Therefore, $\psi(w_n) = \Psi_n(w_n)$, and hence Lemma \ref{l:psin} (c) implies the first identity in \eqref{e:conto}. In order to prove the second identity, we claim that
\begin{equation}
\label{e:lc}
\Psi_k(z_n) = U_k \Psi_{k-1}(z_n) \qquad \mbox{for } 1 \leq k \leq n+1.
\end{equation}
Assume that \eqref{e:lc} holds up to $k-1$, with $1\leq k \leq n+1$. Since $U_j$ maps $(z_{-j-1},z_{-j},z_0,z_j)$ to $(-z_{-j},z_{-j-1},z_j,-z_0)$ and leaves the other coordinates unchanged, the only non-vanishing coordinates of $\Psi_{k-1}(z_n)$ are the one of index $-k$, which has modulus $R+s$, and the one of index $n$ if $k\leq n$, respectively of index $0$ if $k=n+1$, which has value $r$. Therefore,
\[
p_k(\Psi_{k-1}(z_n))=1 \qquad \mbox{and} \qquad q_k(\Psi_{k-1}(z_n)) = \frac{(R+s)^2}{s^2} \geq 1,
\]
so Lemma \ref{l:psin} (a) gives $\Psi_k(z_n) = U_k \Psi_{k-1}(z_n)$. This proves \eqref{e:lc}. The last step $U_{n+1}$ moves the coordinate with value $r$ of index $0$ to index $n+1$ with a minus sign, so 
\begin{equation}
\label{e:last}
\pi_{n+1} \Psi_{n+1}(z_n)=-r, \qquad \pi_j \Psi_{n+1}(z_n) = 0 \quad \forall j\geq 0, \; j\neq n+1,
\end{equation}
and in particular $P_{n+2}(z_n)=0$. By Lemma \ref{l:monotone}, $P_k(z_n)=0$ for every $k\geq n+2$, and Lemma \ref{l:psin} (b) implies that $\Psi_k(z_n)=\Psi_{n+1}(z_n)$ for every $k\geq n+1$. Hence $\psi(z_n)=\Psi_{n+1}(z_n)$, and the second identity in \eqref{e:conto} follows from \eqref{e:last}.
}
\end{rem}

%\bibliographystyle{amsalpha}
%\bibliography{../../biblio/nonlinear}

\newcommand{\etalchar}[1]{$^{#1}$}
\providecommand{\bysame}{\leavevmode\hbox to3em{\hrulefill}\thinspace}
\providecommand{\MR}{\relax\ifhmode\unskip\space\fi MR }
\providecommand{\MRhref}[2]{%
  \href{http://www.ams.org/mathscinet-getitem?mr=#1}{#2}
}
\providecommand{\href}[2]{#2}

\bigskip
\noindent
{\small
Alberto Abbondandolo\\
Ruhr-Universit\"at Bochum, Fakult\"at f\"ur Mathematik\\
Universit\"atsstra{\ss}e 150, 44801 Bochum, Germany\\
\texttt{alberto.abbondandolo@rub.de}

\medskip
\noindent
School of Mathematics, Institute for Advanced Study\\
1 Einstein Drive, Princeton, NJ 08540, USA
}

\end{document}